\documentclass[letterpaper,11pt]{article}

\usepackage{graphicx,epsf,psfrag,amsmath,amssymb}

\newtheorem{theorem}{Theorem}[section]
\newtheorem{lemma}{Lemma}[section]

\newtheorem{definition}{Definition}[section]

\newcommand{\qed}{\hfill\hbox{\rlap{$\sqcap$}$\sqcup$}}
\newenvironment{proof}{\noindent \emph{Proof.\,}}{\qed}

\usepackage{fancyhdr}
\newcommand{\myfootnote}[2]{%
	\thispagestyle{fancy}
	\fancyhf{}
	\renewcommand{\headrulewidth}{0pt}
	\lfoot{\footnotesize\emph{#1}}
	\rfoot{\footnotesize\emph{#2}}
}
\begin{document}

\begin{titlepage}  

\title{On the 1-density of Unit Ball Covering}
\author{
Binhai Zhu
\thanks{Department of Computer Science, Montana State University, Bozeman, MT 59717-3880, USA. Email: {\tt bhz@cs.montana.edu}.}
}

\date{}
\maketitle
\myfootnote{}{December 20, 2007}

\begin{abstract}
Motivated by modern applications like image processing and
wireless sensor networks, we consider a variation of the
Kepler Conjecture. Given any infinite set of unit balls
covering the whole space, we want to know the optimal (lim sup)
density of the volume which is covered by exactly one ball 
(i.e., the maximum such density over all unit ball covers, called the
{\em optimal 1-density} and denoted as $\delta_{d}$, where $d$ is the
dimension of the Euclidean space). We prove that in 2D the optimal
1-density $\delta_{2}=(3\sqrt{3}-\pi)/\pi\approx 0.6539$, which is
achieved through a regular hexagonal covering. In 3D, the problem
is widely open and we present a Dodecahedral Cover Conjecture which
states that the optimal 1-density in 3D, $\delta_3$, is bounded from above
by the 1-density of a unit ball whose Voronoi polyhedron is a regular
dodecahedron of circum-radius one (determined by twelve extra unit balls).
We show numerically that this 1-density $\delta_{3}$({\em dc}) $\approx 0.315$.
\end{abstract}

\end{titlepage}
\newpage

\section{Introduction}

The Kepler conjecture asserts that the density of a packing of unit balls
in three dimensions (3D) is no more than $\pi/\sqrt{18}\approx 0.74048$. This
is known as the oldest (unsolved?) problem in discrete geometry and is
listed as part of the Hilbert's 18th problem. About ten years ago, Hales
published a series of papers concerning the proof of the Kepler conjecture,
using a vast amount of computer verification \cite{Ha98}. While many believe
that the conjecture has been proven, some others believe it has not (as the
referees for {\em Annals of Mathematics} claimed that ``The news from the
referees is bad... They have not been able to certify the correctness of
the proof, and will not be able to certify it in the future" \cite{Sz03}).
In fact, some of Hales' computer programs (total 3GB!) need to re-checked,
which will take another 20 years ({\tt http://en.wikipedia.org/wiki/Kepler\_conjecture}).
\newline

In this paper, we make no attempt in re-proving the Kepler conjecture. In
fact, this research is initiated by looking at the Kepler conjecture from
another perspective: What is the origin of the Kepler conjecture? Application!
It is known that T. Hariot first investigated the Kepler conjecture in 1591
while working for Sir Walter Raleigh, who asked him to determine formulas for
the number of cannonballs regularly packed. Later on Hariot communicated the
problem to Kepler. This problem on unit ball packing is known as the Kepler
conjecture since \cite{Ha98}.
\newline

Looking at the story, it is not difficult for one to realize that
applications do drive theoretical research. In the cannonball packing
problem, obviously cannonballs cannot overlap. But in some modern
applications, the situation could be different. Example 1: in wireless sensor
networks, many homogeneous sensors, whose sensing regions are modeled as
unit disks in 2D, are placed in a region $R$ (we can assume that the diameter
and width of $R$ are much bigger than one). As sensors have limited battery
time, to extend network coverage lifetime (which is the time that the region
is completely covered) one needs to schedule the sensors, i.e., turn some on
and some off. The requirement is that all the time when the network is alive
one must make sure that every point in $R$ is covered by some sensor. The
abstract version of this problem has been shown to be NP-complete \cite{CTLW05},
but the complexity of the geometric version of the problem is still open
(though it seems possible to extend the technique in \cite{CFTZ06} to obtain
a PTAS). In any case, in this wireless sensor network problem clearly the
sensing regions can overlap.
\newline

Another related example is regarding the conversion of continuous-tone images
into binary images for printers, formally known as the {\em digital-halftoning}
technique. Given a continuous-tone image, one can compute
the spatial frequency distribution by applying the Laplacian or Sobel
differential operators. Then we have a grid of points of possibly slightly
different radii (or disks centered at grid points), and the problem is to
compute a subset of disks which maximize the area of the regions covered by
exactly one disk. The complexity of the problem is still open, a factor-5.83
approximation was given in \cite{ABS04}. When the ratio of the radius of
the largest disk over the radius of the smallest disk is a constant, then the
problem admits a PTAS \cite{CFTZ06}.
\newline

\begin{figure}[hbt]
\centerline{\epsffile{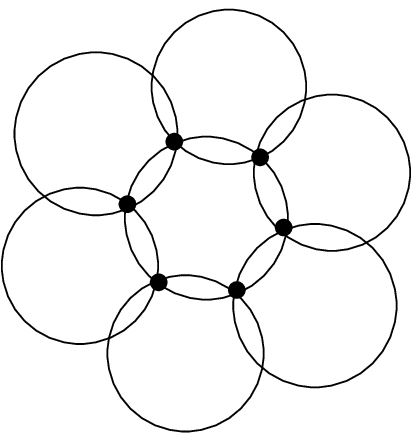}}
\begin{center}{\small Figure 1. A regular hexagonal covering.}
\end{center}
\end{figure}

We feel that the reason why we still could not completely solve these kinds
of geometric problems is that some fundamental question regarding disk
covering, very surprisingly, is still untouched. We cover one such problem
in this paper; namely, the {\em 1-density} of any unit disk cover which covers
all the points in the plane (2D) \footnote{Note that the density for finite
unit disk (ball) covering has long been studied \cite{Bo04}. But this
traditional density is not related to the 1-density we will discuss in this
paper. The traditional density for unit ball covering refers to the maximum
ratio of the volume of some convex region covered by the balls over the volume
of the whole region covered by the balls.}. We prove a theorem which is
corresponding to Thue's theorem on the density of unit disk packing; loosely
speaking, the regular hexagonal covering achieves the optimal 1-density at
0.6539 (Figure 1). (A regular hexagonal covering is the configuration where
the Voronoi polygon of any center of the unit disks is a congruent regular
hexagon. Readers are referred to \cite{PS85} for concepts and properties
regarding Voronoi diagrams.) 
\newline

We then consider the problem in 3D space, i.e., given any unit ball cover
covering all the points in space, decide the corresponding optimal 1-density.
To some extent, this is in fact similar to the Kepler conjecture, which is on
unit ball packing. But the status of this 3D problem is widely open. Our
numerical calculation with computer programs indicates that the 1-density of a
unit ball whose Voronoi polyhedron is a regular dodecahedron with circum-radius
one (determined by twelve extra unit balls) is around 0.315. We conjecture that this
is the upper bound of the optimal 1-density $\delta_3$, which is also in
parallel with the Dodecahedral Conjecture by L.~Fejes~T\'{o}th \cite{To43,HM98}
regarding the Kepler conjecture. We will call this the Dodecahedral Cover
Conjecture henceforth.
\newline

This paper is organized as follows. In Section 2, we present necessary
definitions and a proof for the 2D case. In Section 3, we show how the
1-density in the Dodecahedral Cover Conjecture can be computed (estimated)
empirically. In Section 4, we discuss our conjecture and leave some other
open problems.
 
\section{Result in 2D}

We first make some definitions for any $d$-dimension.
\begin{definition}
A {\em unit ball cover} in $d$-dimensional space is an infinite set of unit
balls which cover every point in the $d$-space.
\end{definition}

\begin{definition}
Given any unit ball cover $C$, let the subset of the balls in $C$ which covers
a larger ball $B$ with a fixed center and a radius $r>1$ be $C(B,r)$. Let
${\cal A}_1(C(B,r))$ be the volume of the parts in $B$ which is covered by
exactly one ball in $C(B,r)$. Let ${\cal A}(B)$ be the volume of $B$. Then,
the {\em 1-density} of $C$, $\delta_{d}(C)$, is defined as
$$\delta_{d}(C)=\lim_{r\rightarrow +\infty}\sup_{B,r}\frac{{\cal A}_1(C(B,r))}{{\cal A}(B)}.$$
Among all possible unit ball covers, the one achieving the maximum 1-density
is denoted as $C^*$ and the corresponding optimal 1-density is denoted as
$\delta_{d}$.
\end{definition}

For the 2D case, a `ball' is usually called a disk and `space' is typically
called a plane. The main theorem for the 2D case is based on the following
lemmas.

\begin{lemma} \label{l1}
Let $P=\{p_1,p_2,p_3,...\}$ be an infinite set of points in the 2D plane
such that $P$ lie on an infinite (or unbounded) number of lines.
Let {\em VD}($P$) be the Voronoi diagram of $P$ such that every point $q$
in the plane is covered by some unit disk centered at $p_j$.
Then the Voronoi polygon of $p_i$, {\em VP}($p_i$), is bounded for all $i$.
\end{lemma}

\begin{proof}
If {\em VP}($p_i$) is unbounded, then $p_i$ is on the convex
hull of $P=\{p_1,p_2,p_3,...\}$. Then some points in {\em VP}($p_i$) is
not covered by any unit disk centered at $p_j$ for any $j$.
A contradiction.
\end{proof}

Let $D=\{D_1,D_2,D_3,...\}$ be an infinite number of disks of a uniform
radius, with each disk $D_i$ centered at $o_i$. Let {\em VD}($D$) be the
Voronoi diagram of the center $o_i$'s for the disks in $D$. For any $o_i$,
let $v_i$ be the Voronoi vertex in the Voronoi polygon {\em VP}($o_i$) which
is the farthest from $o_i$. Define $r_{\max}=\sup_k d(o_k,v_k)$.
We have the following lemma.

\begin{figure}[hbt]
\centerline{\epsffile{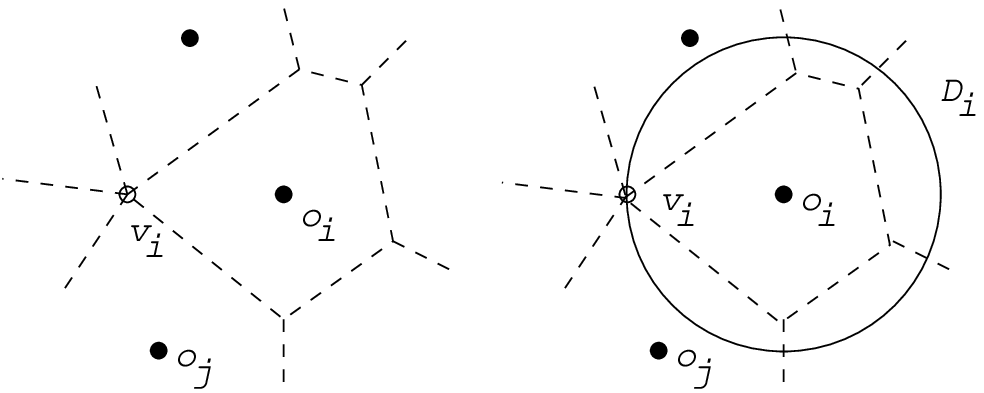}}
\begin{center}{\small Figure 2. Illustration for the proof of Lemma \ref{l2}.}
\end{center}
\end{figure}

\begin{lemma} \label{l2}
$D=\{D_i|i=1,2,3,...\}$ is a unit disk cover if and only if 
$v_i$ is covered by $D_i$, for all $i$.
\end{lemma}

\begin{proof}
We first show the ``if" part.
Without loss of generality, let $r_{\max}=1\geq r_i=d(o_i,v_i)$, for any $i$.
By the definition of $v_i$, a disk $D'_i$ centered at $o_i$ with radius $r_i$
covers {\em VP}($o_i$). As $D_i$ is centered at $o_i$ and with a radius
$r_{\max}=1\geq r_i$, $D_i$ also covers {\em VP}($o_i$) (Figure 2). So every
point in the plane is covered by some unit disk in $D$. $D$ is hence a unit
disk cover.

We now show the ``only if" part.
If $D$ is a unit disk cover, then by definition, every point in the plane
is covered by some disk $D_k$. Following Lemma \ref{l1}, all the Voronoi
polygons {\em VP}($o_i$) in {\em VD}($D$) are bounded, for all $i$.
Without loss of generality, let $r_{\max}=1=d(o_i,v_i)$, for some $i$.
By definition, $v_i$ is the farthest Voronoi vertex in {\em VP}($o_i$)
from $o_i$. Therefore the disk $D_i$ covers $v_i$ (and {\em VP}($o_i$)). 
\end{proof}

\begin{theorem} \label{t1}
$\delta_{2}=(3\sqrt{3}-\pi)/\pi\approx 0.6539$.
\end{theorem}

\begin{proof}
Let $C=\{C_1,C_2,C_3,...\}$ be any unit disk cover, with each disk $C_i$
centered at $o_i$. 
Without loss of generality, let $r_{\max}=1=d(o_i,v_i)$, where $v_i$ is
the farther Voronoi vertex on {\em VP}($o_i$) from $o_i$. By the definition
of a Voronoi vertex, $v_i$ is shared by at least three Voronoi polygons
\footnote{We do not assume that no three centers are collinear, as long as
all the centers are on an infinite number of lines.}. Assume that
one of the other Voronoi polygons containing $v_i$ as a vertex is 
{\em VP}($o_j$), we have $d(o_i,v_i)=d(o_j,v_i)=1$ (notice that
$v_i$ is the farthest Voronoi vertex on {\em VP}($o_j$) from $o_j$, or
we can imagine that $v_j=v_i$).
Following Lemma \ref{l2}, to cover the whole plane we must cover $v_i$
using a unit disk centered at $o_i$. Let $d(o_i,o_j)=2x$.
Denote the sector in $C_i$ bounded by rays $\overrightarrow{o_iv_i}$ and
$\overrightarrow{o_io_j}$ as $\langle v_io_io_j\rangle$ (Figure 3).
Obviously, to maximize the 1-density of sector $\langle v_io_io_j\rangle$,
it is reasonable to assume that no disk $C_k$ intersects the region
in sector $\langle v_io_io_j\rangle$ which is out of the disk $C_j$.
This part of the area in $\langle v_io_io_j\rangle$ which is only covered
by $C_i$, is
$$x\sqrt{1-x^2}-\frac{\arccos x}{2}.$$
By a simple calculation (taking the derivative over $x$ then solving a
simple equation), this function is maximized at
$\frac{3\sqrt{3}-\pi}{12}$ when we set $x=\frac{\sqrt{3}}{2}$, which implies
that the angle between $\overrightarrow{o_iv_i}$ and $\overrightarrow{o_io_j}$
is $\pi/6$ (Figure 3). Correspondingly, the optimal 1-density in the
sector $\langle v_io_io_j\rangle$ is 
$$(\frac{3\sqrt{3}-\pi}{12})/(\frac{\pi}{12})=\frac{3\sqrt{3}-\pi}{\pi},$$
which is roughly 0.6539. By symmetry, the 1-density for $C_i$ is optimal
when $C_j$ and the five other unit disks intersect $C_i$ in a regular hexagonal
setting (e.g., as in Figure 1). Clearly this configuration can be arranged at
all unit disks.
\newline

\begin{figure}[hbt]
\centerline{\epsffile{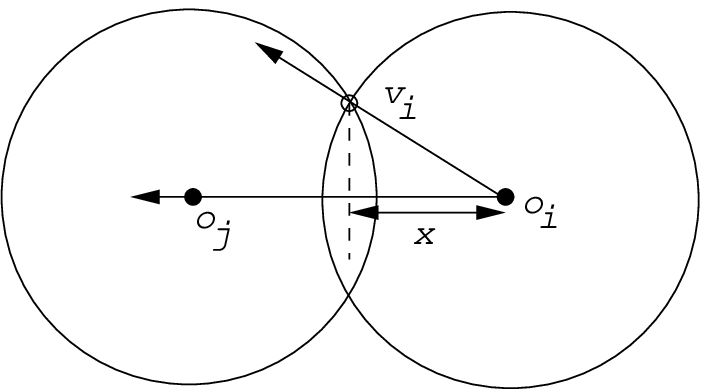}}
\begin{center}{\small Figure 3. Illustration for the proof of Theorem \ref{t1}.}
\end{center}
\end{figure}

Therefore,
$$\delta_{2}=\frac{3\sqrt{3}-\pi}{\pi}\approx 0.6539.$$
\end{proof}

We comment that the above optimal regular hexagonal covering is similar to
Thue's theorem \cite{Thu92,Thu10}, which states that the densest density of unit
disk packing is achieved with a regular hexagonal packing, which admits a
density at $\frac{\pi}{\sqrt{12}}\approx 0.9068$.

\section{Result in 3D}

For the 3D problem, the most relevant application seems to be placing sensors
(or sonars) in the ocean to monitor the temperature change or the appearance of
submarines, etc. It turns out that for the 3D problem, it might be difficult to
apply the same technique as in Theorem \ref{t1} (as the 3D counterpart of the
regular hexagons, i.e., regular dodecahedra, cannot even tile the whole space).
Another reason is the difficulty in computing the volume of a non-convex
body (i.e., it might not even be a polyhedron) analytically. This is very
similar to the Kepler conjecture.
However, different from the Kepler conjecture, the status on
computing/estimating the optimal 1-density $\delta_{3}$ is widely open.
\newline

We conjecture that $\delta_3$ is bounded from above by the 1-density of a unit
ball whose Voronoi polyhedron is a regular dodecahedron with circum-radius one
(determined by twelve extra unit balls). We denote this 1-density as
$\delta_{3}$({\em dc}). As we do not know how to calculate
$\delta_{3}$({\em dc}) analytically (which is related to compute/estimate the
volume of a non-convex body), we first use a random sampling technique
to estimate $\delta_{3}$({\em dc}).
Incidentally, while there have been a lot of theoretical research on estimating
the volume of a convex body, estimating/approximating the volume of a
non-convex body has only been practically investigated recently by the author
of this paper and his colleagues (though without a theoretical guarantee)
\cite{LZZ07}. (The problem of estimating/approximating the volume of a convex
body, with some theoretical guarantee, has been well researched \cite{LV03}.
However, for non-convex geometric body, the problem is still open.) We will
sketch some details of the method (called {\bf Generate-and-Probe}) below
and apply it to obtain an approximate $\delta_{3}$({\em dc}). We comment that
Generate-and-Probe is a very general method, empirically it works even when
no description of the body is given --- as long as there is an oracle which
tells whether a point is inside the body or not.
\newline

Following Figure 4, the unit ball $B_i$ centered at $o_i$ is a
circumscribing ball for
its Voronoi polyhedron {\em VP}($o_i$) --- a regular dodecahedron. Computed by
a computer program, the edge length {\tt a} of {\em VP}($o_i$) is
{\tt a=0.763934}. The distance from $o_i$ to one of the faces of
{\em VP}($o_i$) is {\tt H=0.760071}. Define {\tt R} as the radius of the
circumscribing circle of a face of {\em VP}($o_i$),
which is a pentagon. Then {\tt R=0.649841}. So the ratio $\alpha_1$ of
the volume of {\em VP}($o_i$) to the volume of $B_i$ is
$${\tt \alpha_1=\frac{12*(0.5*a*R*\sin(3.14159*(54.0/180.0))*5)*H*/3}{4*3.14159/3}=0.728762}.$$
To compute the 1-density of $\delta_{3}$({\em dc}), all we need to do
now is to compute the ratio $\alpha_2$ of the volume inside a dodecahedron
which is covered by only one ball to the volume of the dodecahedron. Then
$\delta_3$({\em dc})=$\alpha_1\alpha_2$.
\newline

\begin{figure}[hbt]
\centerline{\epsffile{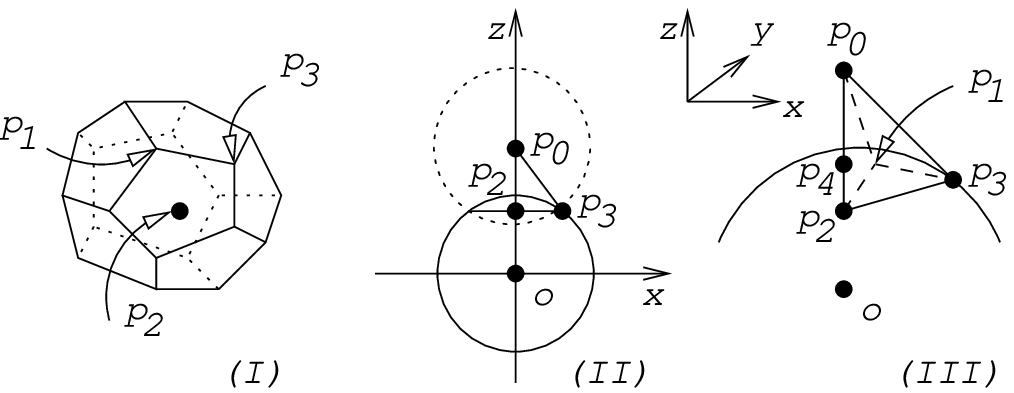}}
\begin{center}{\small Figure 4. Illustration for the computation of $\delta_3$({\em dc}).}
\end{center}
\end{figure}

By symmetry, we consider two unit balls centered
at {\tt o=(0,0,0)} and {\tt $p_0$=(0,0,2H)=(0,0,1.52014)} respectively. The intersection
of the two balls, at height $z=$ {\tt H = 0.760071}, contains a regular
pentagon (face)
of the regular dodecahedron corresponding to the unit ball centered at {\tt $p_0$}.
Let {\tt $p_2$=(0,0,0.760071)} be the center of this pentagon. Again, by
symmetry, we only need to consider a tetrahedron formed by {\tt $p_0,p_2$} and
two other neighboring vertices of the pentagon at height $z=$ {\tt H}. We
choose them as {\tt $p_1$=(0,0.649841,0.760071)},
{\tt $p_3$=(0.618035,0.200812,0.760071)}. Then the volume of the tetrahedron
$T=\diamond(p_0p_1p_2p_3)$ formed by {\tt $p_0,p_1,p_2,p_3$}, {\tt vol(T)},
is analytically computed as
$${\tt vol(T)=0.050877}.$$
The part of $T$, which is out of the ball centered at {\tt $o$=(0,0,0)}, or
{\tt $x^2+y^2+z^2=1$}, is what we need to compute. We denote this part as
$S$ and its volume as {\tt vol(S)}. In Figure 4, we show the bottom-up view
of the regular dodecahedron centered at $p_0$, in (I); the side view along the y-axis
in (II); and the 3D view in (III).
\newline
 
The volume of $S$, which is not even a polyhedron, is hard to compute.
We first use the Generate-and-Probe method of Liu, Zhang and Zhu \cite{LZZ07}.
The method starts by generating $n$ random points in $S$, and then computing
a small probing ball (the radius of the probing ball is selected so that
the probability that the ball contains at least some points in $S$ goes
to 1 when $n\rightarrow +\infty$). Then we use the probing ball to probe
as many parts of $S$ as possible and count the average ratio ${\cal R}$ of
the number of random points enclosed by the probing balls over the volume
of all the probing balls (of course, over a decent number of tries).
Eventually, the volume of $S$ is simply $n/{\cal R}$. We summarize our
empirical results over three cases with $n$=80,000, $n$=140,000 and
$n$=200,000 respectively. The average volume of {\tt vol(S)} is obtained
over 100 tries of the Generate-and-Probe method and the standard deviation
is also given.

\vspace{.1in}

\begin{center}
\begin{tabular}{|c|c|c|c|} \hline
{\tt  n } & {\tt 80,000} & {\tt 140,000} & {\tt 200,000}\\ \hline
avg({\tt vol(S)}) & {\tt 0.0219066} & {\tt 0.0219515} & {\tt 0.0221303}\\ \hline
{\tt $\sigma$(vol(S))} & {\tt 1.66214e-07} & {\tt 6.56446e-08} & {\tt 6.14723e-08}\\ \hline
\end{tabular}
\end{center}
  
\begin{center}
Table 1. The average of {\tt vol(S)} and its standard deviation $\sigma$({\tt vol(S)}), obtained by Generate-and-Probe, over 100 tries.
\end{center}

It remains to verify that the above result by Generate-and-Probe is accurate,
as there is no theoretical guarantee for Generate-and-Probe. We try to verify
this, again empirically. Luckily, the boundary of $S$ can be specified
explicitly. Let the tetrahedron which is obtained by cutting $T$ with a
halfspace $z\geq 1$ be $T'$ (the plane $z=1$ intersects the segment $op_0$
at {\tt $p_4$=(0,0,1)}). Then analytically,
we can compute the volume of $T'$, {\tt vol(T')}. It turns out that
$${\tt vol(T')=0.0163051} < {\tt vol(S)}.$$
The volume of the tetrahedron $\diamond(p_0p_1p_3p_4)$ is
{\tt 0.0348169}, which is an upper bound of {\tt vol(S)}. To obtain a closer
upper bound for {\tt vol(S)}, on the surface of $x^2+y^2+z^2=1$ inside the
tetrahedron $\diamond(p_0p_1p_3p_4)$ (note that {\tt $p_1,p_3,p_4$} are the
three corners of this spherical patch), we generate a set $M$ of $m$ points. We
then compute the 3D convex hull of $M\cup\{p_1,p_3,p_4\}$. We use the upper
hull of {\em CH}($M\cup\{p_1,p_3,p_4\}$) to approximate the patch. The
approximate volume of {\tt vol(S)}, which we still denote as {\tt vol(S)},
is computed by summarizing the volume of all the tetrahedra $\diamond(p_0uvw)$
where $\triangle(uvw)$ is a face on the upper hull of
{\em CH}($M\cup\{p_1,p_3,p_4\}$). The results are summarized in Table 2.

\vspace{.1in}

\begin{center}
\begin{tabular}{|c|c|c|c|} \hline
{\tt  m } & avg({\tt vol(S)}) & {\tt $\sigma$(vol(S))} & $\max$({\tt vol(S)}) \\ \hline
{\tt ~2,000} & {\tt 0.0219467} & {\tt 6.52514e-11} & {\tt 0.0219879} \\ \hline
{\tt ~4,000} & {\tt 0.0219645} & {\tt 1.91297e-11} & {\tt 0.0219956} \\ \hline
{\tt ~6,000} & {\tt 0.0219810} & {\tt 8.92527e-12} & {\tt 0.0219994} \\ \hline
{\tt ~8,000} & {\tt 0.0219843} & {\tt 5.62743e-12} & {\tt 0.0220002} \\ \hline
{\tt 10,000} & {\tt 0.0219873} & {\tt 4.46376e-12} & {\tt 0.0220016} \\ \hline
{\tt 12,000} & {\tt 0.0219914} & {\tt 3.10292e-12} & {\tt 0.0220046} \\ \hline
{\tt 14,000} & {\tt 0.0219931} & {\tt 2.59317e-12} & {\tt 0.0220046} \\ \hline
{\tt 16,000} & {\tt 0.0219947} & {\tt 2.63297e-12} & {\tt 0.0220096} \\ \hline
{\tt 18,000} & {\tt 0.0219966} & {\tt 1.48441e-12} & {\tt 0.0220060} \\ \hline
{\tt 20,000} & {\tt 0.0219970} & {\tt 2.07155e-12} & {\tt 0.0220060} \\ \hline
\end{tabular}
\end{center}
  
\begin{center}
Table 2. The average of {\tt vol(S)}, its standard deviation
$\sigma$({\tt vol(S)}) and its maximum $\max$({\tt vol(S)}), obtained by the 3D Convex Hull method, over 100 tries.
\end{center}

From Table 2, it can also be seen that the maximum of {\tt vol(S)}, over all
of our 1000 tries, has never been greater than {\tt 0.0220096}. Therefore,
with these rigorous computation and verification, it is safe to say that
{\tt vol(S)} $\approx 0.022$ hence $\alpha_2\approx 0.432414$ and
$\delta_3$({\em dc})=$\alpha_1\alpha_2\approx 0.315$. The Convex Hull method
has much smaller deviations, and is probably more accurate in this case, which
is understandable. As we mentioned earlier, the Generate-and-Probe is a more
general (interior-based) approach. Even when no exact description of the body
is given, it still works (at least empirically). We hence have the following
claim.
\newline

\noindent
{\bf Claim.} $\delta_3$({\em dc}) $\approx 0.315$.

\section{Closing Remarks}

The Kepler conjecture on the 3D unit ball packing has been bothering
mathematicians for about 400 years and a long story has been developed
on it. In this paper, based on modern applications in image processing and
wireless sensor networks, we consider the problem of computing the optimal
1-density of a unit ball cover which is a variation of the Kepler conjecture.
We prove a theorem for the 2D case that this optimal 1-density in 2D is
$\delta_{2}=(3\sqrt{3}-\pi)/\pi\approx 0.6539$. The 3D problem is widely open,
though we think that the optimal 1-density is bounded from above by the
1-density of a unit ball whose Voronoi polyhedron is a regular dodecahedral
of circum-radius one (i.e., this Voronoi polyhedron is determined by
twelve extra unit balls). We denote this 1-density as $\delta_3$({\em dc})
and we show numerically that $\delta_3$({\em dc}) $\approx 0.315$. We conclude
this paper with the following formal conjecture.
\newline

\noindent
{\bf The Dodecahedral Cover Conjecture} (for unit ball covering): The optimal
1-density in 3D, $\delta_{3}$, is bounded from above by $\delta_3$({\em dc}) $\approx 0.315$. 
\newline

It is also natural to say that when the dimension $d$ is large, the optimal
1-density in $d$-dimension goes to zero, i.e., $\lim_{d\rightarrow +\infty}\delta_d\rightarrow 0$.

\section*{Acknowledgments}

I thank Sheng Liu for helping running the Generate-and-Probe program to
obtain some of the empirical results. Thanks also to Bin Fu and Marshall Bern
for discussions regarding the Kepler conjecture and the Dodecahedral Conjecture.
Special thanks to Joe O'Rourke for his 3D Convex Hull code in C (on which my
Convex Hull method program is based).

\end{document}